\documentclass[12pt]{article}
\makeatletter
  
  \@addtoreset{equation}{section}
\makeatother
\newtheorem{pp}{Proposition}[section]
\newtheorem{thm}{Theorem}[section]
\newtheorem{lm}{Lemma}[section]

\newtheorem{co}{Corollary}[section]

\usepackage{amsmath,amsfonts,latexsym}
\usepackage{color}
\usepackage{ulem}

\date{ }

\begin{document}

\title{Regularity of Schr\"odinger's functional equation and mean field PDEs for h-path processes
\thanks{To appear in Osaka J. Math., 2019.2010 MSC: Primary 60G30 ; Secondary 93E20}
}

\author{Toshio Mikami\thanks{
This work was  supported by JSPS KAKENHI Grant Numbers JP26400136 and JP16H03948.}\\
\\
}

\maketitle

\begin{abstract}
We show that the solution of Schr\"odinger's functional equation is measurable in space, kernel and marginals.
As an application, we show that the drift vector of the h-path process with given two end point marginals is a measurable function of space, time and marginal at each time.
In particular, we show that the coefficients of mean field PDE systems which  the marginals satisfy 
are measurable function of space, time and marginal.  

\end{abstract}


\section{Introduction}

E. Schr\"odinger considered a probabilistic problem from which he obtained the so-called Schr\"odinger's functional equation
(see section 7 in \cite{S2}  and also \cite{S. B, S1}).
We describe Schr\"odinger's functional equation. 
Let $S$ be a $\sigma$-compact metric space,
let $C(S\times S)$ denote the space of all continuous functions on $S\times S$ with the topology induced by the uniform convergence on every compact subset of $S$
and let ${\cal P}(S)$ denote the space of all Borel probability measures on $S$ with the strong topology.
Fix a positive function $q\in C(S\times S)$.
Schr\"odinger's functional equation can be described as follows.
For $\mu_1, \mu_2\in {\cal P}(S)$,  find a product measure $\nu_1(dx_1)\nu_2(dx_2)$ of 
nonnegative $\sigma$-finite Borel measures on $S$ for which
the following holds:

\begin{equation}\label{1.1}
\begin{cases}
\displaystyle\mu_1(dx_1)=\nu_1(dx_1)\int_{S} q(x_1,x_2)\nu_2(dx_2),&\\
&\\
\displaystyle\mu_2(dx_2)=\nu_2(dx_2)\int_{S} q(x_1,x_2)\nu_1(dx_1)&
\end{cases}
\end{equation}
It is known that (1.1) has the unique solution (see \cite{C, J1} and also \cite{B, Fo}).

\begin{equation}\label{1.2}
u_i(x_i):=\log\biggl(\int_{S} q(x_1,x_2)\nu_j(dx_j)\biggr),\quad i,j =1,2,i\ne j.
\end{equation}
Then $\exp (u_1(x))$ and $\exp (u_2 (x))$ are positive and 
\begin{equation}\label{1.3}
\mu_i (dx)=\exp (u_i (x))\nu_i(dx),\quad i=1,2.
\end{equation}
(1.1)  can be rewritten as follows: for $i,j=1,2$, $i\ne j$,
\begin{equation}\label{1.4}
\exp (u_i (x_i))=\int_{S} q(x_1,x_2)\exp (-u_j(x_j))\mu_j(dx_j),\quad
\mu_i (dx_i)-a.s..
\end{equation}
In particular, Schr\"odinger's problem (1.1) is equivalent to finding a function $u_1(x_1)+u_2(x_2)$ for which 
(\ref{1.4}) holds.
Since $\nu_1(dx_1)\nu_2(dx_2)$ is the unique solution of (1.1),
it is a functional of $\mu_1$, $\mu_2$ and $q$.
Since it is a product measure, $\nu_1$ and $\nu_2$ are also functionals of $\mu_1$, $\mu_2$ and $q$ (see the proof of Corollary 2.1 in section 3):
\begin{equation}\label{1.5}
\nu_i (dx)=\nu_i (dx;q,\mu_1, \mu_2),\quad u_i (x)=u_i (x;q,\mu_1, \mu_2),\quad i=1,2.
\end{equation}
This does not imply the uniqueness of $\nu_1$ and $\nu_2$.
Indeed, for $C>0$, 
$$\nu_1\nu_2=C\nu_1\cdot C^{-1}\nu_2.$$
Let $\{A_n\}_{n\ge 1}$ be an nondecreasing sequence of compact subsets of $S$ such that $S=\cup_{n\ge 1} A_n$. 
$A_1:=S$ when $S$ is compact.
We assume that the following holds so that $\nu_i$, $u_i$, $i=1,2$ are unique:
\begin{equation}\label{1.6}
\quad\nu_1 (A_{n_0 (\mu_1,\mu_2)})=\nu_2(A_{n_0(\mu_1,\mu_2)}),
\end{equation}
where $n_0(\mu_1,\mu_2):=\min \{n\ge 1| \mu_1(A_n)\mu_2(A_n)>0\}.$

Let ${\cal M}(S)$ denote the space of all Radon measures on $S$.
In this paper we denote by a Radon measure a locally finite and inner regular Borel measure.
It is known that a locally finite and $\sigma$-finite Borel measure on a $\sigma$-compact metric space is a Radon measure in our sense (see e.g., p. 901, Prop. 32.3.4 in \cite{G}).

In Theorem 2.1, we show that if $S$ is compact, then
the following are strongly continuous:
$$\nu_i(dx;\cdot,\cdot,\cdot):C(S\times S)\times {\cal P}(S)\times {\cal P}(S) \mapsto{\cal M}(S),$$
$$u_i:C(S\times S)\times {\cal P}(S)\times {\cal P}(S)\mapsto  C(S),$$
and $u_i\in C(S\times C(S\times S)\times {\cal P}(S)\times {\cal P}(S))$.
In Corollary 2.1, we also show that if $S$ is $\sigma$-compact, then
the following are weakly Borel measurable and Borel measurable respectively:
$$\nu_i(dx;\cdot,\cdot,\cdot):C(S\times S)\times {\cal P}(S)\times {\cal P}(S)
\mapsto{\cal M}(S),$$
$$u_i:S\times C(S\times S)\times {\cal P}(S)\times {\cal P}(S)\mapsto  {\bf R}\cup\{\infty\}.$$

As an application of this measurability result,
we show that the coefficients of the mean field PDE system which
the marginal distributions of the h-path process with given two end point marginals
satisfy are measurable functions of space, time and marginal.
To describe the problem more precisely,
we introduce Jamison's result on SDEs for the h-path process with given two end point marginals.
We first describe assumptions and then state Jamison's results.

\noindent
(A1.1) $d\ge 1$ and $\sigma(t,x)=(\sigma^{ij}(t,x))_{i,j=1}^d$, $(t,x)\in [0,1]\times {\bf R}^d$, is a $d\times d$-matrix.
$a(t,x):=\sigma(t,x)\sigma(t,x)^*$, $(t,x)\in [0,1]\times {\bf R}^d$, is uniformly positive definite,
bounded, once continuously differentiable and uniformly H\"older continuous.
$D_x a(t,x)$ is bounded and the first derivatives of $a(t,x)$ are uniformly H\"older continuous
with respect to $x$.

\noindent
(A1.2) $b(t,x):[0,1]\times {\bf R}^d\mapsto {\bf R}^d$ is bounded, continuous and uniformly H\"older continuous
with respect to $x$.
\begin{thm} {\rm (\cite{J2}, p. 330)}
Suppose that (A1.1) and (A1.2) hold.
Then for any $P_0\in {\cal P}({\bf R}^d)$,
 the following SDE has the unique weak solution with a positive continuous transition probability density $p(t,x;s,y)$, $0\le t<s\le 1$, $x,y\in {\bf R}^d$:
\begin{eqnarray}\label{1.7}
dX(t)&=&b(t,X(t))dt+\sigma (t,X(t))dW(t),\quad 0< t<1,\\
PX (0)^{-1}&=&P_0.\nonumber 
\end{eqnarray}
Here $W(t)$ denotes a d-dimensional $\sigma [X(s);0\le s\le t]$-Brownian motion.
Besides, for any $\mu_1, \mu_2\in {\cal P}({\bf R}^d)$, and the solution $\nu_2$ of (1.1)
with $S$ and $q(x_1,x_2)$ respectively replaced  by ${\bf R}^d$ and $p(0,x_1;1,x_2)$,
\begin{equation}\label{1.8}
h(t,x):=\int_{{\bf R}^d} p(t,x;1,x_2)\nu_2 (dx_2)\in C^{1,2} ([0,1)\times {\bf R}^d),
\end{equation}
\begin{equation}\label{1.9}
\left(\frac{\partial}{\partial t}+{\cal A}_t\right) h(t,x)=0, \quad (t,x)\in [0,1)\times {\bf R}^d.
\end{equation}
Here 
$${\cal A}_t:=\frac{1}{2}Trace (a(t,x)D_x^2)
+\langle b(t,x), D_x\rangle .
$$
\end{thm}

\begin{thm}[Markovian reciprocal process] {\rm (\cite{J2}, Theorem 2)}
Suppose that (A1.1) and (A1.2) hold.
Then for any $P_0, P_1\in {\cal P}({\bf R}^d)$ for which $P_1 (dy)\ll dy$, 
there exists  the unique weak solution to the following SDE:
\begin{eqnarray}\label{1.10}
dX(t)&=&\{a(t,X(t))D_x \log h(t,X(t))+b(t,X(t))\}dt\\
&&\qquad +\sigma (t,X(t))dW(t),\quad 0< t< 1,\nonumber \\
PX (t)^{-1}&=&P_t,\quad t=0,1.\nonumber 
\end{eqnarray}
Here, to define $h(t,x)$, we consider (1.1) with $\mu_1, \mu_2, q(x_1,x_2)$ and $S$ respectively replaced by $P_0, P_1, p(0,x_1;1,x_2)$ and 
${\bf R}^d$. $W(t)$ also denotes a d-dimensional $\sigma [X(s);0\le s\le t]$-Brownian motion.
Besides,
\begin{equation}\label{1.11}
PX (t)^{-1}(dx)=\left(\int_{{\bf R}^d} \nu_1 (dx_1)p(0,x_1;t,x)\right) h(t,x)dx,\quad 0\le t\le 1,
\end{equation}
where 
$$\int_{{\bf R}^d} \nu_1 (dx_1)p(0,x_1;0,x)dx:= \nu_1 (dx),$$
$$h(1,x)=\int_{{\bf R}^d} \nu_2 (dx_2)p(1,x;1,x_2):= \frac{\nu_2 (dx)}{dx}.
$$
\end{thm}

\noindent
Remark 1.1.
Replace $S$ by ${\bf R}^d$ in (1.1).
Then the following holds (see (1.2), (1.3), (1.8) and (1.11)): for $x\in {\bf R}^d$,
\begin{equation}\label{1.12}
h(t,x):=\begin{cases}
\exp\{u_1 (x; p(t,\cdot;1,\cdot),PX (t)^{-1},P_1)\},& t\in [0,1),\\
\displaystyle\frac{\nu_2(dx;p(0,\cdot;t,\cdot),P_0,PX (t)^{-1})}{dx}&\\
=\displaystyle \exp\{-u_2 (x; p(0,\cdot;t,\cdot),P_0,PX (t)^{-1})\}\frac{PX (t)^{-1}(dx)}{dx},
&t\in (0,1].
\end{cases}
\end{equation}

As an application of Corollary 2.1 in section 2,  we show that 
\begin{equation}\label{1.13}
U(t,x,P):=\begin{cases}u_1 (x; p(t,\cdot;1,\cdot),P,P_1),\quad t\in [0,1),&\\&\\
\displaystyle\log\left( \frac{\nu_2(dx;p(0,\cdot;t,\cdot),P_0,P)}{dx}\right), \quad t=1&
\end{cases}
\end{equation}
is a Borel measurable function from $[0,1]\times {\bf R}^d\times {\cal P}({\bf R}^d)$ to ${\bf R}$
(see Corollary 2.2).
Theorems 1.1 and 1.2 and 
(\ref{1.12})-(\ref{1.13}) imply that if $P_1(dy)\ll dy$, then $p(t,x)dx:=PX (t)^{-1}(dx)$ satisfies the following mean field PDE system
 (see \cite{Ben, Ben-2, Car, Las} and the references therein for the mean field games and the master equations).
For any $f\in C_b^2 ({\bf R}^d)$ and $t\in (0,1]$,
\begin{eqnarray}\label{1.14}
&&\int_{{\bf R}^d}f(x)p(t,x)dx-\int_{{\bf R}^d}f(x)P_0(dx)\\
&=& \int_0^t ds\int_{{\bf R}^d}({\cal A}_s f(x)+\langle a(s,x)D_x U(s,x,PX (s)^{-1}), Df(x)\rangle )p(s,x)dx,\nonumber
\end{eqnarray}
and for $(t,x)\in (0,1)\times {\bf R}^d$,
\begin{eqnarray}\label{1.142}
0&=&\frac{\partial  U(t,x,PX (t)^{-1}) }{\partial t}+{\cal A}_t U(t,x,PX (t)^{-1})\\
&&\qquad +\frac{1}{2} \langle a(t,x)D_xU(t,x,PX (t)^{-1}), D_xU(t,x,PX (t)^{-1})\rangle,\nonumber
\end{eqnarray}
\begin{eqnarray}
U(1,x,PX (1)^{-1})&=&\log\left (\frac{ \nu_2 (dx;p(0,\cdot;1,\cdot ),P_0, P_1)}{dx}\right).\nonumber
\end{eqnarray}
Here we consider $U(t,x,PX (t)^{-1})$ as a function of $(t,x)$.

Let $\gamma (t;\omega )$ denote a progressively measurable ${\bf R}^d$-valued stochastic process on some filtered probability space and 
consider the following SDE in a weak sense:
\begin{equation}\label{1.15}
dX^\gamma  (t)=\{\gamma (t;\omega )+b(t,X^\gamma (t))\}dt+\sigma (t,X^\gamma  (t)) dW(t),  
\end{equation}
provided it exists  (see e.g. \cite{F}).
Here $W(t)$ denotes a $d$-dimensional Brownian motion defined on the same 
filtered probability space as $\gamma (t;\omega )$.

It is also known that the h-path process with given two end point marginals is the unique minimizer of the following
stochastic optimal control problem
(see \cite{7, 11}, \cite{Leo1}-\cite{28}, \cite{TT}, \cite{Z} and the references therein
for recent progress, especially for stochastic optimal transport).
\begin{thm}{\rm (\cite{7}, \cite{25}, \cite{Z})}
Suppose that (A1.1) and (A1.2) hold.
Then for any $P_0, P_1\in {\cal P}({\bf R}^d)$ for which $P_1 (dy) \ll dy$, 
$\gamma (t;\omega )=a(t,X^\gamma  (t))D_x \log h(t,X^\gamma  (t))$ is the unique minimizer of
the following:
\begin{eqnarray}\label{1.16}
&&V(P_0,P_1)\\
&:=&\inf \biggl\{ E\biggl[\int_0^1 \frac{1}{2}|\sigma (t,X^\gamma (t))^{-1}\gamma (t)|^2dt \biggr]\biggr|
PX ^\gamma (t)^{-1}=P_t ,t=0,1 \biggr\}\nonumber\\
&=&\int_{{\bf R}^d}\log h(1,x)P_1(dx)-\int_{{\bf R}^d}\log h(0,x)P_0(dx),\nonumber
\end{eqnarray}
provided it is finite (see (1.10) for notation).
\end{thm}

\noindent
Remark 1.2.
A sufficient condition for the finiteness of $V(P_0,P_1)$ is given in {\rm \cite{23}} for more general problems.

Schr\"odinger's functional equation (1.1) with $q(x_1,x_2)$ and $S$ respectively 
replaced by $p(0,x_1;1,x_2)$ and ${\bf R}^d$ is equivalent to the Euler equation for  $V(P_0,P_1)$.
We state and prove it for readers' convenience since we could not find any literature
(see Proposition 2.1).

In section 2 we state our main results and prove them in section 3.

\section{Main results}

In this section we state our main results.
We first describe assumptions.

\noindent
(A2.1)
$S$ is a compact metric space.

\noindent
(A2.2)
$q\in C(S\times S;(0,\infty))$.

\noindent
(A2.1)'
$S$ is a $\sigma$-compact metric space.

For a metric space $X$ and $\mu\in {\cal M}(X)$,
\begin{equation}\label{2.1}
||\mu||:=\sup\left\{\int_{X} \phi(x)\mu(dx)\biggl |\phi\in C(X), ||\phi||_\infty\le 1\right\}\in [0,\infty ],
\end{equation}
where for $f\in C(X)$,
\begin{equation}\label{2.2}
||f||_\infty:=\sup_{x\in X}|f(x)|.
\end{equation}

When $S$ is compact, we have the continuity results on $\nu_i$, $u_i$ in (1.5)
(Recall (1.6)). 

\begin{thm}
Suppose that (A2.1) and (A2.2) hold.
Suppose also that $\mu_{i,n}$, $\mu_{i}\in {\cal P}(S)$, $q_n\in C(S\times S;(0,\infty))$, $i=1,2$, $n\ge 1$
and
\begin{equation}\label{2.3}
\lim_{n\to\infty}( ||\mu_{1,n}\mu_{2,n}-\mu_{1}\mu_{2}||+
||q_n-q||_\infty)=0.
\end{equation}
Then
\begin{eqnarray}\label{2.4}
\lim_{n\to\infty}||\nu_1 (\cdot;q_n,\mu_{1,n},\mu_{2,n})\nu_2 (\cdot;q_n,\mu_{1,n},\mu_{2,n})\quad &&\\
\qquad -\nu_1 (\cdot;q,\mu_{1},\mu_{2})\nu_2 (\cdot;q,\mu_{1},\mu_{2})||&=&0,\nonumber\\
\lim_{n\to\infty}\sum_{i=1}^2||u_i (\cdot;q_n,\mu_{1,n},\mu_{2,n})-u_i (\cdot;q,\mu_1, \mu_2)||_\infty&=&0.
\end{eqnarray}
Besides, for $i=1,2$, and $\{x_n\}_{n\ge 1}\subset S$ which converges, as $n\to\infty$, to $x\in S$,
\begin{equation}\label{2.6}
\lim_{n\to\infty}u_i (x_n;q_n,\mu_{1,n},\mu_{2,n})=
u_i (x;q,\mu_1, \mu_2).
\end{equation}
\end{thm}

When $S$ is $\sigma$-compact, we only have the Borel measurability results on $\nu_i$, $u_i$ in (1.5). 

\begin{co}
Suppose that (A2.1)' and (A2.2) hold.
Then the following are Borel measurable: for $i=1,2$,
$$\int_S f(x)\nu_i(dx;\cdot,\cdot, \cdot):C(S\times S)\times {\cal P}(S)\times {\cal P}(S)
\mapsto {\bf R},\quad f\in C_0 (S),$$
$$u_i:S\times C(S\times S)\times {\cal P}(S)\times {\cal P}(S)
\mapsto  {\bf R}\cup\{\infty\}.$$
\end{co}

As an application of Corollary 2.1, we obtain the following.

\begin{co}
Suppose that (A1.1) and (A1.2) hold.
Then $U(t,x,P)$ in (1.13) is a Borel measurable function from $[0,1]\times {\bf R}^d\times {\cal P}({\bf R}^d)$ to ${\bf R}$.
In particular, 
(\ref{1.14})-(\ref{1.142}) hold.
\end{co}

For $P_0\in {\cal P}({\bf R}^d)$ and Borel measurable $f: {\bf R}^d\mapsto {\bf R}$,
\begin{equation}\label{2.8}
V_{P_0}^*(f):=\sup\biggl\{\int_{{\bf R}^d}f(x)P(dx)-V(P_0,P):P\in {\cal P}({\bf R}^d)\biggr\}
\end{equation}
(see (\ref{1.16}) for notation).
Then since $P\mapsto V(P_0,P)$ is convex, lower semicontinuous and $\not\equiv \infty$, 
for $P\in {\cal P}({\bf R}^d)$,
\begin{equation}\label{2.7}
V(P_0,P)=\sup\biggl\{\int_{{\bf R}^d}f(x)P(dx)-V_{P_0}^*(f):f\in C_b({\bf R}^d)\biggr\}\in [0,\infty ]
\end{equation}
(see \cite{21-1,22, 25, TT} and the references therein).
The following gives the variational meaning to Schr\"odinger's functional equation.

\begin{pp}
Suppose that (A1.1) and (A1.2) hold.
Then for any $P_0, P_1\in {\cal P}({\bf R}^d)$ for which $P_1 (dy) \ll dy$
and for which $V(P_0,P_1)$ is finite,
Schr\"odinger's functional equation (1.1) with $\mu_1$, $\mu_2$ and $q(x_1,x_2)$ 
respectively replaced by $P_0$, $P_1$ and $p(0,x_1;1,x_2)$ is equivalent to the following:
\begin{equation}\label{2.8}
P_1 (dy)=\frac{\delta V_{P_0}^*(\log h(1,\cdot))}{\delta f} (dy).
\end{equation}
Here  $\displaystyle\frac{\delta V_{P_0}^*(f)}{\delta f} $ denotes the G\^ateaux derivative of $V_{P_0}^*(f)$.
\end{pp}


\section{Proof of main results}
 
 In this section we state and prove lemmas and prove our main results.
 \begin{eqnarray}\label{3.1}
m_q&:=&\inf \{q(x_1,x_2)|x_1,x_2\in S\},\\
M_q&:=&\sup \{q(x_1,x_2)|x_1,x_2\in S\}.\nonumber
\end{eqnarray}
The following two lemmas are proved in \cite{B}.
\begin{lm}{\rm (\cite{B}, p. 194)}
Suppose that (A2.1) and (A2.2) hold.
Then, for any $\mu_1, \mu_2\in {\cal P}(S)$,
there exists a unique pair of nonnegative finite measures $\nu_{1}$, $\nu_{2}$ on $S$ for which (\ref{1.1}) and the following holds:
\begin{equation}\label{3.2}
\frac{1}{\sqrt{M_q}}\le\nu_{1}(S)=\nu_{2}(S)\le \frac{1}{\sqrt{m_q}},
\end{equation}
\begin{equation}\label{3.3}
\frac{m_q}{\sqrt{M_q}}\le  \exp( u_i (x))
\le \frac{M_q}{\sqrt{m_q}}, \quad x\in S, i=1,2
\end{equation}
(see (1.2) for notation).
\end{lm}

\begin{lm}{\rm (\cite{B}, section 7)}
Suppose that (A2.1) and (A2.2) hold.
Then, there exists a function $c(a,b)$ which is nonincreasing in $a$ and nondecreasing in $b$
such that
for any $\mu_i, \tilde\mu_i\in {\cal P}(S)$, $i=1,2$,
\begin{equation}\label{3.4}
||\nu_1\nu_{2}-\tilde\nu_{1}\tilde\nu_{2}||\le c(m_q,M_q)
||\mu_{1}\mu_{2}-\tilde\mu_{1}\tilde\mu_{2}||^{\frac{1}{2}}.
\end{equation}
Here $\tilde\nu_{i}(dx):=\nu_{i }(dx;q,\tilde\mu_1,\tilde\mu_2)$, $i=1, 2$
(see (1.5) and (2.1) for notation).
\end{lm}

The following lemma can be proved by Lemma 3.1.
\begin{lm}
Suppose that (A2.1) and (A2.2) hold and that $q_n\in C(S\times S;(0,\infty))$, $n\ge 1$
and 
\begin{equation}\label{3.5}
\lim_{n\to\infty}||q_n-q||_\infty =0
\end{equation}
(see (2.2) for notation).
Then, for any $\mu_i\in {\cal P}(S)$, $i=1,2$,
\begin{equation}\label{3.6}
\lim_{n\to \infty} ||\nu_{n,1}\nu_{n,2}-\nu_{1}\nu_{2}||=0,
\end{equation}
where $\nu_{n,i }(dx):=\nu_{i }(dx;q_n,\mu_1,\mu_2)$
(see (1.5) and (2.1) for notation).
\end{lm}
(Proof)
$u_{n,i}(x):=u_{i}(x;q_n ,\mu_1,\mu_2)$.
Then, from (\ref{1.2})-(\ref{1.3}), 
\begin{eqnarray}\label{3.7}
u_{n,i} (x_i)&=&\log \left(\int_{S} q_n (x_1,x_2)\nu_{n,j }(dx_j)\right),\quad i,j=1,2, i\ne j,\\
\nu_{n,i }(dx)&=&\exp(-u_{n,i } (x))\mu_{i }(dx), \quad i=1,2.\nonumber
\end{eqnarray}
For $i=1,2$, $\displaystyle \left\{\frac{\nu_{n,i}(dx)}{\nu_{n,i}(S)}\right\}_{n\ge 1}$ is a tight family of 
probability measures and $\{\nu_{n,i}(S)\}_{n\ge 1}$ is bounded from above and below by (\ref{3.2}). 
In particular, there exist $\{s(n)\}_{n\ge 1}$ and  a finite measure $\overline \nu_{i}$
such that $\nu_{s(n),i }$ weakly converges, as $n\to\infty$, to $\overline\nu_{i }$.
From construction, (3.2) with $\nu_i$ replaced by $\overline\nu_{i }$ also holds.
\begin{equation}\label{3.8}
\overline u_{i} (x_i):=\log\left(\int_{S} q (x_1,x_2)\overline \nu_{j}(dx_j)\right),\quad i,j=1,2, i\ne j.
\end{equation}
Then for $i=1,2$,
\begin{equation}\label{3.9}
\overline\nu_{i }(dx)=\exp(-\overline u_{i}(x))\mu_{i}(dx).
\end{equation}
Indeed, from (3.7),
\begin{eqnarray}
&&\nu_{s(n),i}(dx)-\exp(-\overline u_{i} (x))\mu_{i}(dx)\nonumber\\
&=&(\exp (-u_{s(n),i} (x))-\exp(-\overline u_{i} (x)))\mu_{i}(dx).\nonumber
\end{eqnarray}
For $i,j=1,2, i\ne j$ and $x_i\in S$,
\begin{eqnarray}
&&|\exp (u_{s(n),i} (x_i))-\exp(\overline u_{i} (x_i))|\nonumber\\
&\le&\left|\int_S (q_{s(n)} (x_1,x_2)-q (x_1,x_2))\nu_{s(n),j }(dx_j)\right|\nonumber\\
&&\qquad +\left|\int_S q (x_1,x_2)(\nu_{s(n),j }(dx_j)-\overline\nu_{j }(dx_j))\right|\nonumber\\
&\le&||q_{s(n)} -q||_\infty\times \nu_{s(n),j}(S)\nonumber\\
&&\qquad +\left|\int_S q (x_1,x_2)(\nu_{s(n),j }(dx_j)-\overline \nu_{j }(dx_j))\right|
\to0, \quad n\to\infty,\nonumber
\end{eqnarray}
from (3.2) and (3.5).
From (\ref{3.3}),
$$\exp (-u_{s(n),i } (x_i))\le 
\frac{\sqrt{M_{q_{s(n)} }}}{m_{q_{s(n)} }}\to \frac{\sqrt{M_{q }}}{m_{q }},  \quad n\to\infty, i=1,2.
$$
In particular, the bounded convergence theorem implies that  (3.9) is true.

From (3.8)-(3.9),
\begin{eqnarray}\label{3.10}
\mu_{i}(dx_i)&=&\exp(-\overline u_{i}(x_i))\mu_{i}(dx_i)\exp(\overline u_{i}(x_i))\\
&=&\overline\nu_{i }(dx_i)\int_S q(x_1,x_2)\overline\nu_{j }(dx_j), \quad i,j=1,2, i\ne j.
\nonumber
 \end{eqnarray}
The uniqueness of the solution to (\ref{1.1}) implies that 
\begin{equation}\label{3.11}
\overline\nu_{i }(dx)=\nu_{i }(dx),\quad i=1,2
\end{equation}
since (3.2) hold for both of $\overline\nu_{i }$ and $\nu_{i }$.
Since the above method applies for any subsequence of $\{q_n\}_{n\ge 1}$,
the discussion in (3.9) implies that the following holds:
\begin{equation}\label{3.12}
\lim_{n\to\infty}
||\nu_{n,i}-\nu_{i }||=0, \quad i=1,2.
\end{equation}
(\ref{3.2}) and (\ref{3.12}) completes the proof.$\Box$

We prove Theorem 2.1 by Lemmas 3.1-3.3.

\noindent
(Proof of Theorem 2.1)

\noindent
Lemmas 3.2 and 3.3 imply (2.4).
We prove (2.5).
Without loss of generality, we only have to consider the case when $i=1$.
For sufficiently large $n\ge 1$, 
\begin{eqnarray}\label{3.13}
&&||u_1(\cdot ;q_n,\mu_{1,n},\mu_{2,n})-u_1 (\cdot;q,\mu_1, \mu_2)||_\infty\\
&\le &-\log\biggl\{ 1-\frac{\sqrt {M_q}}{m_q}\biggl ( \frac{||q_n-q||_\infty}{\sqrt {m_{q_n}}}\nonumber\\
&&\qquad + ||q||_\infty\cdot ||\nu_2 (\cdot;q_n,\mu_{1,n},\mu_{2,n})-\nu_2  (\cdot;q,\mu_{1},\mu_{2})||\biggl)\biggr\}\nonumber\\
&\to &0, \quad n\to\infty.\nonumber
\end{eqnarray}
We prove (3.13).
For $x\in S$,
\begin{eqnarray}
&&u_1(x;q_n,\mu_{1,n},\mu_{2,n})-u_1 (x;q,\mu_1, \mu_2)\nonumber\\
&=&\log \left(1+\frac{\int_S q_n(x,x_2)\nu_2 (dx_2;q_n,\mu_{1,n},\mu_{2,n})-\int_S q(x,x_2)\nu_2  (dx_2;q,\mu_{1},\mu_{2})}{
\int_S q(x,x_2)\nu_2  (dx_2;q,\mu_{1},\mu_{2})}\right),\nonumber
\end{eqnarray}
\begin{eqnarray}
&&\left|\frac{\int_S q_n(x,x_2)\nu_2 (dx_2;q_n,\mu_{1,n},\mu_{2,n})-\int_S q(x,x_2)\nu_2  (dx_2;q,\mu_{1},\mu_{2})}{
\int_S q(x,x_2)\nu_2  (dx_2;q,\mu_{1},\mu_{2})}\right|\nonumber\\
&\le & \frac{1}{\int_S q(x,x_2)\nu_2  (dx_2;q,\mu_{1},\mu_{2})}
\biggl\{\left|\int_S (q_n(x,x_2)-q(x,x_2))\nu_2 (dx_2;q_n,\mu_{1,n},\mu_{2,n})\right|\nonumber\\
&&\qquad +\left|\int_S q(x,x_2)(\nu_2 (dx_2;q_n,\mu_{1,n},\mu_{2,n})-\nu_2  (dx_2;q,\mu_{1},\mu_{2}))\right|\bigg\}\nonumber\\
&\le &\frac{\sqrt {M_q}}{m_q}\biggr(\frac{||q_n-q||_\infty}{\sqrt {m_{q_n}}}+ ||q||_\infty\cdot ||\nu_2 (\cdot;q_n,\mu_{1,n},\mu_{2,n})-\nu_2  (\cdot;q,\mu_{1},\mu_{2})||\biggr)\nonumber
\end{eqnarray}
from(3.2)-(3.3).
The following also holds:
$$\lim_{n\to\infty}||\nu_2 (\cdot;q_n,\mu_{1,n},\mu_{2,n})-\nu_2  (\cdot;q,\mu_{1},\mu_{2})||=0.
$$
Indeed, for $f\in C_b (S)$ for which $||f||_\infty\le 1$, from (2.4) and (3.2),
\begin{eqnarray}
&&\int_S f(x)(\nu_2 (dx;q_n,\mu_{1,n},\mu_{2,n})-\nu_2  (dx;q,\mu_{1},\mu_{2}))\nonumber\\
&=& \frac{1}{\nu_1 (S;q_n,\mu_{1,n},\mu_{2,n})}
\int_S f(x)\{ (\nu_1 (S;q_n,\mu_{1,n},\mu_{2,n})\nu_2 (dx;q_n,\mu_{1,n},\mu_{2,n})\nonumber\\
&&\qquad\qquad -\nu_1  (S;q,\mu_{1},\mu_{2})\nu_2  (dx;q,\mu_{1},\mu_{2}))\nonumber\\
&&\qquad +\left(\nu_1 (S;q,\mu_{1},\mu_{2})-\nu_1 (S;q_n,\mu_{1,n},\mu_{2,n})\right)\nu_2  (dx;q,\mu_{1},\mu_{2})\}\nonumber\\
&\le&  \sqrt{M_{q_n}}\biggl\{||\nu_1 (\cdot;q_n,\mu_{1,n},\mu_{2,n})\nu_2 (\cdot;q_n,\mu_{1,n},\mu_{2,n})\nonumber\\
&&\qquad\qquad -\nu_1  (\cdot;q,\mu_{1},\mu_{2})\nu_2  (\cdot;q,\mu_{1},\mu_{2})||\nonumber\\
&&\qquad +\biggl|
\sqrt{\nu_1 (S;q,\mu_{1},\mu_{2})\nu_2 (S;q,\mu_{1},\mu_{2})}\nonumber\\
&&\qquad-
\sqrt{\nu_1 (S;q_n,\mu_{1,n},\mu_{2,n})\nu_2 (S;q_n,\mu_{1,n},\mu_{2,n})}\biggl|\frac{1}{\sqrt {m_q}}\biggr\}\nonumber\\
&\to& 0, \quad n\to\infty.\nonumber
\end{eqnarray}
The following completes the proof of (3.13):
$$
\log (1-|a|)\le \log (1+a)\le \log (1+|a|)\le -\log  (1-|a|),\quad |a|<1.
$$
We prove (2.6).
From (2.5), we only have to prove the following: 
for $i=1,2$, and $\{x_n\}_{n\ge 1}\subset S$ which converges to $x\in S$ as $n\to\infty$, 
\begin{equation}\label{3.14}
\lim_{n\to\infty}u_i (x_n;q,\mu_{1},\mu_{2})=
u_i (x;q,\mu_1, \mu_2).
\end{equation}
This can be proved by the bounded convergence theorem.$\Box$ 

For $\mu_1, \mu_2\in {\cal P}(S)$,
\begin{equation}\label{3.15}
\mu_{i|n} (E):=\frac{\mu_{i} (E\cap A_n)}{\mu_{i} (A_n)}, \quad E\in {\cal B}(S), \quad n\ge 
n_0(\mu_1, \mu_2),  i=1,2,
\end{equation}
where ${\cal B}(S)$ denotes the Borel $\sigma$-field of $S$ (see (\ref{1.6}) for notation).
When we replace $X$ and $S$ by $A_n$ in (2.1)-(2.2) and (3.1), 
we use notations $||\cdot||_n$, $||\cdot||_{\infty,n}$, $m_{q,n}$ and $M_{q,n}$
instead of $||\cdot||$, $||\cdot||_{\infty}$, $m_{q}$ and $M_{q}$ respectively.
We use a similar convention when it is not confusing.

We introduce and prove two lemmas to prove Corollary 2.1.

\begin{lm}
Suppose that (A2.1)' and (A2.2) hold.
Then, for any $\mu_1, \mu_2\in {\cal P}(S)$ and any  $k\ge n_0 (\mu_1, \mu_2)$,
there exists a unique pair of nonnegative finite measures $\nu_{1|k}$, $\nu_{2|k}$ on $A_k$ for which  Lemma 3.1 with $S$, $m_{q}$, $M_{q}$, $\mu_{i}$, $\nu_{i}$, $u_i$, $i=1,2$
replaced by $A_k$, $m_{q,k}$, $M_{q,k}$, $\mu_{i|k}$, $\nu_{i|k}$, $u_{i|k}$, $i=1,2$ respectively holds.
Suppose, in addition, that $\mu_{i,n}\in {\cal P}(S)$, $q_n\in C(S\times S;(0,\infty))$, $i=1,2$, $n\ge 1$
and
$$\lim_{n\to\infty}(||\mu_{1,n}\mu_{2,n}-\mu_{1}\mu_{2}||_k+||q_n-q||_{\infty, k})=0.$$
Then (2.4)-(2.6) hold even if $\nu_{i}$, $\mu_{i,n}$, $\mu_{i}$, $||\cdot ||$,
$u_i$, $||\cdot||_\infty$ and $S$ is replaced by 
$\nu_{i|k}$, $\mu_{i,n|k}$, $\mu_{i|k}$, $||\cdot ||_k$, $u_{i|k}$, 
$||\cdot||_{\infty,k}$ and $A_k$ respectively. 
\end{lm}
(Proof)
Theorem 2.1 and the following completes the proof:
\begin{equation}\label{3.16}
||\mu_{1,n|k}\mu_{2,n|k}-\mu_{1|k}\mu_{2|k}||_k
\le 2\frac{||\mu_{1,n}\mu_{2,n}-\mu_{1}\mu_{2}||_k}{\mu_{1}(A_k)\mu_{2}(A_k)}.
\end{equation}
(\ref{3.16}) is true, since
\begin{eqnarray}
&&\mu_{1,n|k}(dx_1)\mu_{2,n|k}(dx_2)- \mu_{1|k}(dx_1) \mu_{2|k}(dx_2)\nonumber\\
&=&\frac{1}{\mu_{1}(A_k)\mu_{2}(A_k)}\{(\mu_{1,n}(dx_1)\mu_{2,n}(dx_2)-\mu_{1}(dx_1)\mu_{2}(dx_2))\nonumber\\
&&\qquad+\frac{\mu_{1}(A_k)\mu_{2}(A_k)-\mu_{1,n}(A_k)\mu_{2,n}(A_k)}
{\mu_{1,n}(A_k)\mu_{2,n}(A_k)}\mu_{1,n}(dx_1)\mu_{2,n}(dx_2)\}.\Box\nonumber
\end{eqnarray}

For any $\mu_1, \mu_2\in {\cal P}(S)$ and any  $n\ge n_0 (\mu_1, \mu_2)$,
\begin{equation}\label{3.17}
\mu^{(n)}(dx_1dx_2):=q(x_1,x_2)1_{A_n \times A_n }(x_1,x_2)\nu_{1|n}(dx_1)
\nu_{2|n}(dx_2).
\end{equation}
The following is known.
\begin{lm}{\rm(\cite{J1}, Theorem 3.2)}
Suppose that (A2.1)' and (A2.2) hold.
Then for any $\mu_i\in {\cal P}(S)$, $i=1,2$,
there exists a unique solution $\nu_{1}(dx_1)\nu_{2}(dx_2)$ to (1.1) and 
$\mu^{(n)}(dx_1dx_2)$ weakly converges, as $n\to\infty$, to $\mu (dx_1dx_2):=q(x_1,x_2)\nu_{1}(dx_1)\nu_{2}(dx_2)$.
\end{lm}

By Lemmas 3.4 and 3.5, we prove Corollary 2.1.

\noindent
(Proof of Corollary 2.1)
Without loss of generality, we only have to prove the case when $i=1$.
From Lemma 3.4, for any $n\ge 1$,
$\int_{S\times S}f(x_1,x_2)\mu^{(n)}(dx_1dx_2)$ is continuous in $(f, q, \mu_1,\mu_2)$
on the open set
$$C_b(S\times S)\times C(S\times S) \times \{(\mu_1,\mu_2)\in {\cal P}(S)\times {\cal P}(S)|
n_0(\mu_1,\mu_2)\le n\}$$
(see (\ref{1.6})).
Notice  that $n_0(\mu_1,\mu_2)\le n$ if and only if $\mu_1(A_n)\mu_2(A_n)>0$.
From Lemma 3.5, $\int_{S\times S}f(x_1,x_2)\mu (dx_1dx_2)$ is measurable in $(f, q, \mu_1,\mu_2)$.
The following implies the first part of Corollary 2.1: for $f\in C_o (S)$ and $x\in\mathbb{R}$,
\begin{align*}
&\biggl\{(q,\mu_1,\mu_2)\in C(S\times S)\times {\cal P}(S)\times {\cal P}(S)\biggl|
\int_{S}f(x_1)\nu_1(dx_1;q,\mu_1,\mu_2)<x\biggl\}\\
=&\bigcup_{k=1}^\infty
\biggl\{(q,\mu_1,\mu_2)\in C(S\times S)\times {\cal P}(S)\times {\cal P}(S)\biggl|
n_0(\mu_1,\mu_2)=k,\notag\\
&\qquad \frac{\int_{S\times A_k }f(x_1)q(x_1,x_2)^{-1}\mu(dx_1dx_2;q,\mu_1,\mu_2)}
{\sqrt{\int_{A_k\times A_k }q(x_1,x_2)^{-1}\mu(dx_1dx_2;q,\mu_1,\mu_2)}}
<x\biggr\}.\notag
\end{align*}
For any $x\in S$ and $\phi_n \in C_o (S)$ for which $0\le \phi_n\le 1$ and  $\phi_n (y)=1, y\in A_n$,
$\int_{S}\phi_n (x_2)q(x,x_2)\nu_2(dx_2)$
is measurable in $(q,\mu_1, \mu_2)$  in the same way as above and is continuous in $x$.
In particular, it is measurable in $(x, q,\mu_1, \mu_2)$ and so is the following: by Fatou's lemma,
$$\exp (u_{1} (x))=\lim_{n\to\infty}\int_{S}\phi_n (x_2)q(x,x_2)\nu_2(dx_2).\Box$$

Corollary 2.1 immediately implies Corollary 2.2.

\noindent
(Proof of Corollary 2.2)
Since $p(t,\cdot;1,\cdot)$ is continuous on $[0,1)$ from Theorem 1.1,
$$(t,x,P, P_1)\mapsto (x,p(t,\cdot;1,\cdot), P,P_1)$$
is continuous on $[0,1)\times  {\bf R}^d\times {\cal P}({\bf R}^d)\times {\cal P}({\bf R}^d)$, which 
implies the measurability of $U(t,x,P)$.
It is easy to see that (1.14) - (1.15) hold.$\Box$

We prove Proposition 2.1.

\noindent
(Proof of Proposition 2.1)
\begin{equation}\label{3.18}
\frac{\delta V_{P_0}^* (\log h(1,\cdot))}{\delta f}(dy)=h(1,y)dy\int_{{\bf R}^d} p(0,x;1,y)
\frac{P_0 (dx)}{h(0,x)}.
\end{equation}
Indeed, for any $\psi \in C_b({\bf R}^d)$ and $\varepsilon\in {\bf R}$, instead of $P_1$,
consider Schr\"odinger's problem (1.1) with $S, \mu_1, \mu_2$ and $q(x_1,x_2)$ 
respectively replaced by ${\bf R}^d, P_0, \mu^{\varepsilon\psi}$ and $q(0,x_1;1,x_2)$,
where for $f \in C_b({\bf R}^d)$
$$\mu^{f}(dy):=h(1,y)\exp (f(y))dy \int_{{\bf R}^d} P_0 (dx) \frac{p(0,x;1,y)}{\int_{{\bf R}^d} p(0,x;1,z)h(1,z)\exp (f(z))dz}$$
(see (1.1)-(1.4)).
Then, from Theorem 1.3 and (2.8) (see e.g. \cite{7, Z} and also \cite{25}), 
\begin{eqnarray}
&&V_{P_0}^* (\log h(1,\cdot)+\varepsilon \psi)\nonumber\\
&=&\int_{{\bf R}^d} \log \left(\int_{{\bf R}^d} p(0,x;1,y)h(1,y)\exp (\varepsilon \psi (y))dy \right)P_0 (dx)\nonumber
\end{eqnarray}
(see (1.8)).
This implies (3.18).
From (1.8), 
\begin{equation}\label{3.19}
P_0 (dx)=\left(\int_{{\bf R}^d}h(1,y)dy\hbox{}p(0,x;1,y)\right)\frac{P_0 (dx)}{h(0,x)}.
\end{equation}
(3.18) and (3.19) completes the proof.$\Box$



Department of Mathematics, Tsuda University,

2-1-1 Tsuda-machi, Kodaira, Tokyo 187-8577, Japan


\begin{thebibliography}{00}


\bibitem{Ben}
A. Bensoussan, J. Frehse, P. Yam:
Mean Field Games and Mean Field Type Control Theory,
SpringerBriefs in Mathematics, Springer, Berlin-Heidelberg-New York, 2013.

\bibitem{Ben-2}
A. Bensoussan, J. Frehse, P. Yam:
\textit{On the interpretation of the Master Equation},
Stochastic Process. Appl. \textbf{127} (2017), 2093-2137.

\bibitem{S. B}
S. Bernstein: 
\textit{Sur les liaisons entre les grandeurs al\'etoires},
Verh. des intern. Mathematikerkongr. Zurich 1932, Band \textbf{1}, 288-309.

\bibitem{B}
A. Beurling:
\textit{An Automorphism of Product Measures},
Ann. of Math. \textbf{72} (1960), 189-200.

\bibitem{Car}
P. Cardaliaguety:
\textit{Notes on Mean Field Games
(from P.-L. Lions' lectures at College de France)}, January 15, 2012.

\bibitem{C}
Y. Chen, T. Georgiou, M. Pavon:
\textit{Entropic and displacement interpolation: a computational approach using the Hilbert metric},
SIAM J. Appl. Math. \textbf{76} (2016), 2375-2396.

\bibitem{7} 
P. Dai Pra:
\textit{A stochastic control approach to reciprocal diffusion processes}
Appl. Math. Optim. \textbf{23} (1991), 313-329.

\bibitem{F}
W.~H. Fleming,  H.~M. Soner:
Controlled Markov Processes and Viscosity Solutions, 2nd ed.,
Springer, Berlin-Heidelberg-New York, 2006.

\bibitem{11} 
H. F\"ollmer:
\textit{Random fields and diffusion processes},
in: Hennequin, P.~L.(ed) 
\'Ecole d'\'Et\'e de Probabilit\'es de Saint-Flour XV--XVII, 1985--87,  Lecture Notes in Math.\textbf{1362}, 101--203, 
 Springer, Berlin-Heidelberg-New York, 2006.
  
\bibitem{Fo} 
R. Fortet:
\textit{R\'esolution d'un Syst\`eme d'\'equations de M. Schroedinger},
J. Math. Pures Appl. \textbf{IX} (1940), 83-105.

\bibitem{G} 
D. J. H. Garling:
A Course in Mathematical Analysis: Volume 3, Complex Analysis, Measure and Integration,
Cambridge University Press, Cambridge, 2014.

\bibitem{J1}
B. Jamison:
\textit{Reciprocal Processes},
Z. Wahr. verw. Gebiete \textbf{30} (1974), 65--86.

\bibitem{J2} 
B. Jamison:
\textit{The Markov process of Schr\"odinger},
Z. Wahr. Verw. Gebiete \textbf{32} (1975), 323--331.

\bibitem{Las} 
 J-M. Lasry,  P. L. Lions:
\textit{Mean field games},
Japan J. Math. \textbf{2} (2007), 229-260.

\bibitem{Leo1}
C. L\'eonard :  
\textit{From the Schr\"odinger problem to the Monge-Kantorovich problem},
J. Funct. Anal. \textbf{262} (2012), 1879--1920.

\bibitem{Leo2} 
C. L\'eonard :  
\textit{A survey of the Schr\"odinger problem and some of its connections with optimal transport},
Special Issue on Optimal Transport and Applications,
Discrete Contin. Dyn. Syst. \textbf{34} (2014), 1533--1574.


\bibitem{19}  
T. Mikami:
\textit{Monge's problem with a quadratic cost by the zero-noise limit of $h$-path  processes},
Probab. Theory Related Fields \textbf{129} (2004), 245--260.


\bibitem{21-1}
T. Mikami:
\textit{Marginal problem for semimartingales via duality},
in: Giga, Y. , Ishii, K., Koike, S. et al. (eds)
International Conference for the 25th Anniversary of Viscosity Solutions, 
Gakuto International Series. Mathematical Sciences and Applications \textbf{30}, 133--152,
Gakkotosho, Tokyo, 2008.

\bibitem{22} 
T. Mikami:
\textit{Optimal transportation problem as stochastic mechanics},
in Selected Papers on Probability and Statistics, Amer. Math. Soc. Transl. Ser. 2,  \textbf{227}, 75--94, 
Amer. Math. Soc., Providence, RI., 2009.

\bibitem{23}
T. Mikami:
\textit{Two end points marginal problem by stochastic optimal transportation},
SIAM J. Control Optim. \textbf{53} (2015), 2449-2461.

\bibitem{25}
T. Mikami, M. Thieullen: 
\textit{Duality theorem for the stochastic optimal control problem},
Stochastic Process. Appl. \textbf{116} (2006), 1815--1835.
 
\bibitem{28}
 L. R\"uschendorf,  W. Thomsen:
\textit{Note on the Schr\"odinger equation and $I$-projections},
Statist. Probab. Lett. \textbf{17}  (1993), 369--375.

\bibitem{S1}
E. Schr\"odinger: 
\textit{Ueber die Umkehrung der Naturgesetze},
Sitz. Ber. der Preuss. Akad. Wissen., Berlin, Phys. Math., p. 144, 1931.

\bibitem{S2}
E. Schr\"odinger: 
\textit{Th\'eorie relativiste de l'electron et l'interpr\'etation de la m\'ecanique quantique},
Ann. Inst. H. Poincar\'e \textbf{2} (1932), 269-310. 

\bibitem{TT}
X. Tan, N. Touzi:
\textit{Optimal transportation under controlled stochastic dynamics},
Ann. Probab. \textbf{41} (2013), 3201--3240.


\bibitem{Z}
J.~C. Zambrini: 
\textit{Variational processes},
in: Albeverio, S., etal. (eds.)
Stochastic processes in classical and quantum systems (Ascona, 1985), Lecture Notes in Phys. \textbf{262}, 517--529, Springer, Berlin-Heidelberg-New York, 1986.

\end{thebibliography}
\end{document}